\documentclass{article}[11pt]

\usepackage[ruled,vlined]{algorithm2e}

\usepackage{graphicx}
\usepackage{amsmath,amsfonts,amssymb,latexsym}
\usepackage{enumerate}
\usepackage{setspace}
\usepackage{color}
\usepackage[cp1250]{inputenc}
\usepackage{tikz}

\newtheorem{thm}{Theorem}[section]
\newtheorem{prop}[thm]{Proposition}
\newtheorem{ob}[thm]{Observation}

\newtheorem{cor}[thm]{Corollary}

\newtheorem{lemma}[thm]{Lemma}

\newtheorem{prob}{Problem}
\newtheorem{quest}{Question}
\newcommand{\qed}{$\Box$}
\newcommand{\proof}{\noindent\textbf{Proof. }}

\newcommand{\cH}{{\cal H}}
\newcommand{\cN}{{\cal N}}
\newcommand{\cT}{{\cal T}}
\newcommand{\cO}{{\cal O}}
\newcommand{\diam}{{\rm diam}}

\newcommand{\smallqed}{{\tiny ($\Box$)}}

\newcommand{\3}{ \vspace{0.3cm} }
\newcommand{\2}{ \vspace{0.2cm} }

\newenvironment{unnumbered}[1]{\trivlist \item [\hskip \labelsep {\bf
#1}]\ignorespaces\it}{\endtrivlist}

\def\vertex(#1){\put(#1){\circle*{1.8}}}
\def\lab(#1)#2{\put(#1){\makebox(0,0)[c]{#2}}}

\newcommand{\rhogr}{\rho_{\rm gr}}
\newcommand{\taugr}{\tau_{\rm gr}}
\newcommand{\gt}{\gamma_t}
\newcommand{\gr}{\gamma_{\rm gr}}
\newcommand{\grt}{\gamma_{\rm gr}^t}
\newcommand{\modo}{{\rm mod}}
\newcommand{\gtd}{\gamma_{\rm tg}}

\begin{document}

\title{Total Dominating Sequences in Graphs}

\author{
Bo\v{s}tjan Bre\v{s}ar$^{a,b}$
\and
Michael A. Henning$^{c}$
\and
Douglas F. Rall$^{d}$\footnote{Correspondence to: Department of Mathematics, Furman University,
3300 Poinsett Highway, Greenville, SC, 29613, USA.}
}

\date{\today}

\maketitle

\begin{center}
$^a$ Faculty of Natural Sciences and Mathematics, University of Maribor, Slovenia (Bostjan.Bresar@um.si)\\
$^b$ Institute of Mathematics, Physics and Mechanics, Ljubljana, Slovenia\\
$^c$ Department of Pure and Applied Mathematics, University of Johannesburg, South Africa (mahenning@uj.ac.za)\\
$^d$ Department of Mathematics, Furman University, Greenville, SC, USA (doug.rall@furman.edu)\\

\end{center}

\begin{abstract}
A vertex in a graph totally dominates another vertex if they are adjacent. A sequence of vertices in a graph $G$ is called a total dominating sequence if every vertex $v$ in the sequence totally dominates at least one vertex that was not totally dominated by any vertex that precedes $v$ in the sequence, and at the end all vertices of $G$ are totally dominated.
While the length of a shortest such sequence is the total domination number of $G$, in this paper we investigate total dominating sequences of maximum length, which we call the Grundy total domination number, $\grt(G)$, of $G$.
We provide a characterization of the graphs $G$ for which $\grt(G)=|V(G)|$ and of those for which $\grt(G)=2$.
We show that if $T$ is a nontrivial tree of order~$n$ with no vertex with two or more leaf-neighbors, then  $\grt(T) \ge \frac{2}{3}(n+1)$, and characterize the extremal trees.
We also prove that for $k \ge 3$, if $G$ is a connected $k$-regular graph of order~$n$ different from $K_{k,k}$, then $\grt(G) \ge  (n + \lceil \frac{k}{2} \rceil - 2)/(k-1)$ if $G$ is not bipartite and $\grt(G) \ge (n + 2\lceil \frac{k}{2} \rceil - 4)/(k-1)$ if $G$ is bipartite.
The Grundy total domination number is proven to be bounded from above by two times the Grundy domination number, while the former invariant can be arbitrarily smaller than the latter.
Finally, a natural connection with edge covering sequences in hypergraphs is established, which in particular yields the NP-completeness of the decision version of the Grundy total domination
number.
\end{abstract}



\newpage
\section{Introduction}

The concept of edge covering sequences was introduced in~\cite{bgm-2014} to shed more light on the possible procedures of determining the edge cover number of a hypergraph (\emph{edge cover number} is the cardinality of a smallest set of (hyper)edges in a hypergraph whose union equals the set of its vertices). Of particular interest is the maximum length of a sequence, in which one only uses the most basic greedy condition that each edge must contain a vertex that is not contained in the edges that precede it, and is called the Grundy covering number of a hypergraph. (The name arises from the Grundy coloring number, which is the maximum number of colors that are used in a greedy coloring algorithm. The concept of Grundy colorings was introduced back in the 1970's~\cite{cs-79} and has been investigated in many papers.) In~\cite{bgm-2014} the main focus was on dominating sequences (of vertices) in graphs, which can be viewed precisely as edge covering sequences of the hypergraph of closed neighborhoods of the graph. The longest possible dominating sequences were determined in several classes of graphs (e.g. trees, split graphs, cographs), while it was shown that this problem is NP-complete, even when restricted to chordal graphs~\cite{bgm-2014}.

In this paper we introduce and investigate total dominating sequences in graphs, which arise from the hypergraph of open neighborhoods of a graph. Total domination is one of the classical concepts in graph theory, having numerous applications and connections with other parameters. It was recently surveyed in the monograph~\cite{MHAYbookTD}.
The \emph{total domination number}, $\gt(G)$, of a graph $G$ with no isolated vertices is the smallest cardinality of a set of vertices $S$ such that every vertex of $G$ has a neighbor in $S$. (If the condition only requires that vertices from $V(G)\setminus S$ have a neighbor in $S$, then the resulting invariant is the \emph{domination number} $\gamma(G)$ of $G$.) It is well-known that for every graph $G$ with no isolated vertices we have $\gamma(G)\le \gt(G)\le 2\gamma(G)$.
One of the central problems in this area is to determine good upper bounds for the total domination number of a graph in terms of its order. Cockayne, Dawes, and Hedetniemi~\cite{CoDaHe80} showed that if $G$ is connected of order~$n \ge 3$, then $\gt(G) \le \frac{2}{3}n$. Several authors~\cite{Alfewy,ChMc,Tuza} showed that if $G$ is a graph of order~$n$ with minimum degree at least~$3$, then $\gt(G) \le \frac{1}{2}n$. Thomasse and Yeo~\cite{ThYe07} showed that if $G$ is a graph of order~$n$ with minimum degree at least~$4$, then $\gt(G) \le \frac{3}{7}n$.

We now introduce our main invariant, which is defined for all graphs $G$ without isolated vertices.
Let $S=(v_1,\ldots,v_k)$ be a sequence of distinct vertices of $G$. The corresponding set $\{v_1,\ldots,v_k\}$ of vertices from the sequence $S$ will be denoted by $\widehat{S}$. The sequence $S$  is a {\em legal (open neighborhood) sequence} if
\begin{equation}
\label{e:total}
N(v_i) \setminus \bigcup_{j=1}^{i-1}N(v_j) \ne\emptyset.
\end{equation}
holds for every $i\in\{2,\ldots,k\}$.
If, in addition, $\widehat{S}$ is a total dominating set of $G$, then we call $S$ a \emph{total dominating sequence} of $G$. If $S$ is a legal sequence, then we will say that $v_i$ \emph{footprints} the vertices from $N(v_i) \setminus \cup_{j=1}^{i-1}N(v_j)$, and that $v_i$ is the \emph{footprinter} of every vertex $u\in N(v_i) \setminus \cup_{j=1}^{i-1}N(v_j)$. That is, $v_i$ footprints vertex $u$ if $v_i$ totally dominates $u$, and $u$ is not totally dominated by any of the vertices that precede $v_i$ in the sequence. Thus the function $f_S \colon V(G)\to \widehat{S}$ that maps each vertex to its footprinter is well defined.
Clearly the length $k$ of a total dominating sequence $S$ is bounded from below by the total domination number, $\gamma_t(G)$, of $G$. On the other hand, the maximum length of a total dominating sequence in $G$ will be called the \emph{Grundy total domination number} of $G$ and will be denoted by $\grt(G)$. The corresponding sequence will be called a \emph{Grundy total dominating sequence} of $G$.

The paper is organized as follows. In the next section we fix the notation and state some preliminary
results and observations. In particular we prove an upper bound for the Grundy total domination number
in terms of the order and minimum degree of a graph, and a lower bound in terms of the order and maximum degree. Section~\ref{sec:chain} considers two total domination chains that arise from some invariants related to the Grundy total domination number, notably the total domination number, the game total domination number, and the upper total domination number. In Section~\ref{sec:large} we characterize two extremal families of graphs, that is, the graphs whose Grundy total domination number is equal to 2, and the graphs whose Grundy total domination number is equal to their order. While the former are exactly complete multipartite graphs, the latter family can only be described in a more involved fashion, which in the class of trees reduces to exactly the trees having a perfect matching; this result is established in Section~\ref{sec:trees}. This section also contains the proof of the lower bound $\grt(T) \ge \frac{2}{3}(n+1)$, where $T$ is an arbitrary tree, together with the characterization of the trees attaining this bound. Section~\ref{sec:regular} contains our most involved result, which is the lower bound for the Grundy total domination number of regular graphs, when complete bipartite graphs are excluded. In Section~\ref{sec:relations} the bounds between the Grundy total domination number and the Grundy domination number are discussed, while Section~\ref{sec:hyper} connects the new concept with edge covering sequences of hypergraphs. As a result of these connections, we first  establish the existence of total dominating sequences in $G$ of arbitrary length between $\gamma_t(G)$
and $\grt(G)$, and then we prove the NP-completeness of the corresponding Grundy total domination problem. We conclude in the last section with some open problems that arise throughout the paper.

\section{Notation and Preliminary Results}
\label{sec:prelim}

For notation and graph theory terminology, we in general follow~\cite{MHAYbookTD}. We assume throughout the remainder of the paper that all graphs considered are \emph{without isolated vertices}. The \emph{degree} of a vertex $v$ in $G$, denoted $d_G(v)$, is the number of neighbors, $|N_G(v)|$, of $v$ in $G$. The minimum and maximum degree among all the vertices of $G$ are denoted by $\delta(G)$ and $\Delta(G)$, respectively. A \emph{leaf} is a vertex of degree~$1$, while its neighbor is a \emph{support vertex}. A \emph{strong support vertex} is a vertex with at least two leaf-neighbors. The subgraph induced by a set $S$ of vertices of $G$ is denoted by $G[S]$. A \emph{non-trivial graph} is a graph on at least two vertices.

A \emph{cycle} on $n$ vertices is denoted by $C_n$ and a \emph{path} on $n$ vertices by $P_n$.
A \emph{star} is a tree $K_{1,n}$ for some $n \ge 1$.
A \emph{complete $k$-partite graph} is a graph that can be partitioned into $k$ independent sets, so that every pair of vertices from two different independent sets are adjacent. A \emph{complete multipartite graph} is a graph that is complete $k$-partite for some $k$. In particular, complete bipartite and complete graphs are in the family of complete multipartite graphs.

Two distinct vertices $u$ and $v$ of a graph $G$ are \emph{open twins} if $N(u)=N(v)$. A graph is \emph{open twin-free} if it has no open twins. We remark that a tree is open twin-free if and only if it has no strong support vertex.

A \emph{total dominating set} of a graph $G$ with no isolated vertex is a set $S$ of vertices of $G$ such that every vertex is adjacent to a vertex in $S$; that is, every vertex has a neighbor in $S$. If we only require that every vertex outside $S$ has a neighbor in $S$, then $S$ is called a \emph{dominating set} of $G$. The \emph{upper total domination number}, $\Gamma_t(G)$, of $G$ is the maximum cardinality of a minimal total dominating set in $G$.

Given a subset $X$ of vertices in a graph $G$, a legal sequence $S$ of $G$ is a \emph{total dominating sequence} of $X$ if $\widehat{S}$ totally dominates the set $X$ and each vertex of the sequence $S$ footprints a vertex of $X$ not footprinted by any vertex preceding it in $S$. In particular, if $X = V(G)$, then $S$ is a total dominating sequence of $G$.

For a matching $M$ in a graph $G$ a vertex incident to an edge of $M$ is called \emph{strong} if its degree is~$1$ in the subgraph $G[V(M)]$.  The matching $M$ is called a \emph{strong matching} (also called an \emph{induced matching} in the literature) if every vertex in $V(M)$ is strong. The number of edges in a maximum strong matching of $G$ is the \emph{strong matching number}, $\nu_s(G)$, of $G$. The strong matching number is studied, for example, in~\cite{JRS-2014,KMM-2012}. As defined in~\cite{GH-2005}, $M$ is a \emph{semistrong} matching if every edge in $M$ has a strong vertex. The number of edges in a maximum semistrong matching of $G$ is the \emph{semistrong matching number}, $\nu_{ss}(G)$, of $G$.

We are now in a position to present some preliminary results and observations on the Grundy total domination number of a graph. Recall that all graphs in this paper have no isolated vertex. Let $G$ be a graph, and let $H$ be an induced subgraph of $G$ that contains no isolated vertex. Every Grundy total dominating sequence in $H$ is either a total dominating sequence of $G$ or can be extended to a total dominating sequence of $G$, implying that $\grt(G) \ge \grt(H)$. This implies the following result.

\begin{ob}
For every graph $G$, $\grt(G) \ge \max \, \{ \grt(H) \}$, where the maximum is taken over all induced subgraphs $H$ of $G$ with no isolated vertex.
\label{induced}
\end{ob}

If $M$ is a maximum strong matching in a graph $G$, then the subgraph, $H = G[V(M)]$, of $G$ induced by the edges of $M$ is isomorphic to $\nu_s(G)$ disjoint copies of $K_2$, implying by Observation~\ref{induced} that $\grt(G) \ge \grt(H) = 2\nu_s(G)$. Thus, the Grundy total domination number of a graph is at least twice its strong matching number.

\begin{ob}
For every graph $G$, $\grt(G) \ge 2\nu_s(G)$.
\label{indmatch}
\end{ob}

We present next the following general lower bound on the Grundy total domination number of a graph in terms of its order and maximum degree.

\begin{prop}
If $G$ is a graph of order~$n$ with maximum degree~$\Delta(G) = \Delta$, then $\grt(G) \ge \frac{n}{\Delta}$. Further, if $G$ is connected and $\grt(G) = \frac{n}{\Delta}$, then $G = K_{\Delta,\Delta}$.
\label{prop:Delta}
\end{prop}
\proof The lower bound follows immediately from the observation that $\grt(G) \ge \gt(G)$ and the well-known observation (see,~\cite{MHAYbookTD}) that $\gt(G) \ge n/\Delta$.
Suppose that $G$ is connected and $\grt(G) = n/\Delta$.
Let $S$ be an arbitrary total dominating sequence of $G$ and let $|S| = k$. The set $\widehat{S}$ is a total dominating set of $G$. Consequently, $\grt(G) \ge k \ge \gt(G) \ge n/\Delta$. This implies that every total dominating sequence is a Grundy total dominating sequence. As $S$ is a Grundy total dominating sequence and $k = n/\Delta$, every vertex $v$ in $\widehat{S}$ footprints exactly~$\Delta$ vertices. It follows that $G$ is $\Delta$-regular. We show that $G = K_{\Delta,\Delta}$. Suppose, to the contrary, that $G \ne K_{\Delta,\Delta}$. Let $v_1$ be an arbitrary vertex of $G$. Since $G$ is a connected $\Delta$-regular graph, there exists a vertex $v_2$ in $G$ different from $v_1$ with the property that $v_2$ has a neighbor in $N(v_1)$ and a neighbor not in $N(v_1)$. But then there exists a total dominating sequence of $G$ starting with the vertices $v_1$ and $v_2$ as its first two vertices such that $v_2$ footprints strictly less than~$\Delta$ vertices, a contradiction. Therefore, if $\grt(G) = n/\Delta$, then $G = K_{\Delta,\Delta}$.~\qed

\medskip
The following general upper bound on the Grundy total domination number of a graph is in terms of its order and minimum degree.

\begin{prop}
\label{p:deg}
If $G$ is a graph of order~$n$, then $\grt(G) \le n - \delta(G) + 1$.
\end{prop}
\proof
Let $S=(s_1,\ldots, s_k)$ be a Grundy total dominating sequence of $G$.
Let $u$ be a vertex footprinted in the last step, that is,  $u\in f_S^{-1}(s_k)$.
Since $u$ is not totally dominated before the last step, we have $N(u)\cap \{s_1,\ldots,s_{k-1}\}=\emptyset$,
and so
\[
|\{s_1,\ldots,s_{k-1}\}| = k-1 \le n-d(u).
\]
Thus, $\grt(G)=k \le n-\delta(G)+1$.
\qed

\medskip
The upper bound from Proposition~\ref{p:deg} is clearly achieved by complete graphs and by the graph $2K_3+e$. Note that $|V(G)|-\delta(G)+1$ can be at most $|V(G)|$, which is achieved when $\delta(G)=1$. Graphs $G$ with $\grt(G) = |V(G)|$ will be studied in Section~\ref{sec:large}.

\section{Total Domination Chains}
\label{sec:chain}

If $S$ is a sequence of vertices in a graph $G$ such that $\widehat{S}$ is a minimal total dominating set in $G$ of maximum cardinality $\Gamma_t(G)$, then $S$ is a total dominating sequence of $G$ since each vertex in $S$ footprints, among other vertices, the vertices that it uniquely totally dominates in the set $\widehat{S}$. This implies that $\Gamma_t(G) \le \grt(G)$. This gives rise to the following total domination chain.

\begin{ob}
For every graph $G$, $\gamma_t(G) \le \Gamma_t(G) \le \grt(G)$.
\label{chain}
\end{ob}

A natural problem is to characterize the connected graphs for which we have equality throughout the inequality chain given in the statement of Observation~\ref{chain}; that is, graphs $G$ for which $\gt(G) = \grt(G)$. In Section~\ref{sec:large}, we characterize graphs $G$ with $\gt(G) = \grt(G) = 2$. These are shown (in Theorem~\ref{thm:multi}) to be precisely the complete multipartite graphs. The following result shows that there is no graph $G$ satisfying $\gt(G) = \grt(G) = 3$.

\begin{thm}
If $G$ is a graph with $\gt(G) = 3$, then $\grt(G) > 3$.
\label{thm:Tdom3}
\end{thm}
\proof Let $G$ be a graph with $\gt(G) = 3$. Let $S = \{a,b,c\}$ be a minimum total dominating set of $G$. Since $G[S]$ contains no isolated vertex, either $G[S] = P_3$ or $G[S] = K_3$. If $G[S] = P_3$, then renaming vertices if necessary, we may assume that $G[S]$ is the path $abc$. If $G[S] = K_3$, then
$G[S]$ is the $3$-cycle $abca$. In both cases, by the minimality of the set $S$, there is a vertex $a'$ outside $S$ that is adjacent to $a$ but to no other vertex of $S$.
Suppose that $G[S] = P_3$. Let $H$ be the subgraph of $G$ induced by $\{a',a,b,c\}$. Then, $H$ is isomorphic to $P_4$, and by Observation~\ref{induced} and the observation that $\grt(P_4) = 4$, $\grt(G) \ge \grt(H) = 4$. Hence we may assume that every minimum total dominating set in $G$ induces a $K_3$, for otherwise $\grt(G) > 3$ as desired. By assumption, the set $S' = \{a',a,b\}$ which induces a path $P_3$ is not a total dominating set in $G$. Let $c'$ be a vertex not totally dominated by $S'$. Since $S$ is a total dominating set of $G$, this implies that $c'$ is a vertex outside $S$ that is adjacent to $c$ but to no other vertex of $S$. But then $a'acc'$ is an induced $P_4$ in $G$, implying once again that $\grt(G) \ge \grt(P_4) = 4$.~\qed

\begin{figure}[htb]
\begin{center}
\begin{tikzpicture}[scale=.8,style=thick,x=1cm,y=1cm]
\def\vr{2.5pt} 
\path (0,0) coordinate (v1);
\path (.8,0) coordinate (v2);
\path (2,0) coordinate (v3);
\path (2.8,0) coordinate (v4);
\path (4,0) coordinate (v5);
\path (4.8,0) coordinate (v6);
\path (0,2) coordinate (w1);
\path (.8,2) coordinate (w2);
\path (2,2) coordinate (w3);
\path (2.8,2) coordinate (w4);
\path (4,2) coordinate (w5);
\path (4.8,2) coordinate (w6);
\path (6.5,0) coordinate (a1);
\path (7.3,0) coordinate (a2);
\path (8.1,0) coordinate (a3);
\path (8.9,0) coordinate (a4);
\path (6.5,2) coordinate (b1);
\path (7.3,2) coordinate (b2);
\path (8.1,2) coordinate (b3);
\path (8.9,2) coordinate (b4);
%
\draw (v1) -- (w3); \draw (v1) -- (w4); \draw (v1) -- (w5); \draw (v1) -- (w6);
\draw (v2) -- (w3); \draw (v2) -- (w4); \draw (v2) -- (w5); \draw (v2) -- (w6);
\draw (v3) -- (w1); \draw (v3) -- (w2); \draw (v3) -- (w5); \draw (v3) -- (w6);
\draw (v4) -- (w1); \draw (v4) -- (w2); \draw (v4) -- (w5); \draw (v4) -- (w6);
\draw (v5) -- (w1); \draw (v5) -- (w2); \draw (v5) -- (w3); \draw (v5) -- (w4);
\draw (v6) -- (w1); \draw (v6) -- (w2); \draw (v6) -- (w3); \draw (v6) -- (w4);
\draw (a1) -- (b2); \draw (a1) -- (b3); \draw (a1) -- (b4);
\draw (a2) -- (b1); \draw (a2) -- (b3); \draw (a2) -- (b4);
\draw (a3) -- (b1); \draw (a3) -- (b2); \draw (a3) -- (b4);
\draw (a4) -- (b1); \draw (a4) -- (b2); \draw (a4) -- (b3);
\draw (v1) [fill=white] circle (\vr); \draw (v2) [fill=white] circle (\vr); \draw (v3) [fill=white] circle (\vr);
\draw (v4) [fill=white] circle (\vr); \draw (v5) [fill=white] circle (\vr); \draw (v6) [fill=white] circle (\vr);
\draw (w1) [fill=white] circle (\vr); \draw (w2) [fill=white] circle (\vr); \draw (w3) [fill=white] circle (\vr);
\draw (w4) [fill=white] circle (\vr); \draw (w5) [fill=white] circle (\vr); \draw (w6) [fill=white] circle (\vr);
\draw (a1) [fill=white] circle (\vr); \draw (a2) [fill=white] circle (\vr);
\draw (a3) [fill=white] circle (\vr); \draw (a4) [fill=white] circle (\vr);
\draw (b1) [fill=white] circle (\vr); \draw (b2) [fill=white] circle (\vr);
\draw (b3) [fill=white] circle (\vr); \draw (b4) [fill=white] circle (\vr);
\draw (2.4,-0.7) node {(a) \hskip .1cm $G_3$};
\draw (7.7,-0.7) node {(b) \hskip .1cm $G_4$};
\end{tikzpicture}
\end{center}
\vskip -0.7cm
\caption{The graphs $G_3$ and $G_4$.} \label{f:prism}
\end{figure}
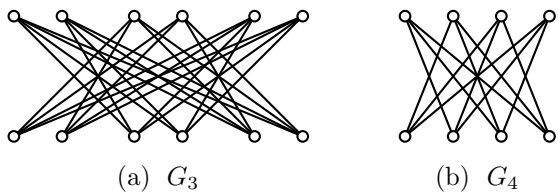

\medskip
Infinite families $\{\mathcal{G}_m\}_{m\ge 3}$ of connected graphs with both total domination number and Grundy total domination number equal to~$4$ can be constructed as follows.  Let $m$ be an integer such that $m \ge 3$.  For each $i$ with $1 \le i \le m$,
let $X_i$ and $Y_i$ be nonempty sets of vertices such that the sets $X_1,\ldots,X_m,Y_1,\ldots,Y_m$
are pairwise disjoint.  Let $X=\cup_{i=1}^m X_i$ and $Y=\cup_{i=1}^m Y_i$.  A bipartite graph $G$ with
$V(G) = X \cup Y$ is obtained by adding the  edge  $xy$ if and only if $x\in X_i$ and $y \in Y_j$
for some $i$ and $j$ such that $i\not= j$.   It is easy to see that $\gt(G)=4$.  Furthermore, every total dominating
sequence  $S$ of $G$ satisfies  $\widehat{S}=\{a,b,c,d\}$ where $a\in X_i$, $b\in X_j$,
$c\in Y_r$, $d \in Y_s$, for any choice of $\{i,j,r,s\}$ such that $i\not= j$ and $r \not= s$.  In fact, any
permutation of such a set of four vertices is a total dominating sequence of $G$.  Hence, $\gt(G)=4=\grt(G)$.
Define $\mathcal{G}_m$ to be the class of all such graphs constructed in this way. We note that the $6$-cycle is the smallest graph in the family~$\mathcal{G}_3$. The graphs $G_3 \in \mathcal{G}_3$ and $G_4\in \mathcal{G}_4$ shown in Figure~\ref{f:prism}(a)
and~\ref{f:prism}(b), are examples of this construction.
We state our observation formally as follows.

\begin{ob}
There are infinitely many connected graphs $G$ with $\grt(G) = \gt(G) = 4$.
\label{ob:gt4}
\end{ob}

The domination game in graphs was introduced in~\cite{brklra-2010} and extensively studied afterwards (see, for example,~\cite{dkr-15,bill-2013}). The total version of the domination game was investigated in~\cite{hekara-2013,hekara-2015}. This game is played on a graph $G$ by two players, named Dominator and Staller. They alternately take turns choosing vertices of $G$ such that each chosen vertex totally dominates at least one vertex not totally dominated by the vertices previously chosen. Dominator's goal is to totally dominate the graph as fast as possible, and Staller wishes to delay the process as long as possible. The game total domination number, $\gtd(G)$, of $G$ is the number of vertices chosen when Dominator starts the game and both players play optimally. Every sequence of vertices generated by Dominator and Staller in the total domination game is a total dominating sequence, implying the following result.

\begin{ob}
For every graph $G$ with no isolated vertex, $\gt(G) \le \gtd(G) \le \grt(G)$.
\label{chain_game}
\end{ob}

We remark that the difference between the game total domination number and the Grundy total domination number can be arbitrarily large. For example, for $k \ge 2$, if $G$ is the graph of order~$n = 2k+1$ obtained from $k$ disjoint copies of $K_3$ by identifying one vertex from each copy of $K_3$ into a common vertex (of degree~$2k$), then $\gtd(G) = 2$ and $\grt(G) = n-1$.

\section{Graphs with large or small Grundy total domination number}
\label{sec:large}

In this section we provide a characterization of the graphs $G$ for which $\grt(G)=|V(G)|$ and of those for which $\grt(G)=2$. The latter value of $\grt(G)$ is the least possible, since $\grt(G)\ge \gamma_t(G)\ge 2$. We begin with a lemma that will be used in characterizing the graphs of order~$n$ that also have  Grundy total domination number~$n$.

\begin{lemma}  \label{lem:reciprocal}
Let $G$ be a graph of order $n$ such that $\grt(G)=n$ and let $S=(v_1,\ldots,v_n)$ be a
Grundy total dominating sequence of $G$.  If $x$ and $y$ are any two vertices of $G$ such that
$x$ footprints $y$ with respect to $S$, then $y$ also footprints $x$ with respect to $S$.
\end{lemma}
\proof  Since $S=(v_1,\ldots,v_n)$ is a total dominating sequence of $G$, it follows that each vertex of $G$ footprints exactly one vertex with respect to $S$.  This means that the footprinter function $f_S:V(G)\to \widehat{S}$ is injective.  Let $\overrightarrow{G}$ be the directed graph that has vertex set $V(G)$ and that has $\{(f_S(u),u) \mid u\in V(G)\}$ as its set of directed edges.  In $\overrightarrow{G}$ each vertex has in-degree 1 and out-degree 1.  Consequently, $\overrightarrow{G}$ is the disjoint union of directed cycles.  Let  $x_1 \ldots x_k x_1$ be any one of these directed cycles. We may assume without loss of generality that $x_1$ is the first of the vertices in $\{x_1,\ldots,x_k\}$
to appear in $S$.   In particular, $x_1$ footprints $x_2$ and $x_k$ since $x_1$ is adjacent to both $x_2$ and  $x_k$. That is, $f_S(x_2)=x_1=f_S(x_k)$. We conclude that $k=2$.  Hence, $\overrightarrow{G}$ is the disjoint union of directed cycles of order 2.  Thus, if a vertex $x$ footprints a vertex $y$, then $y$ also footprints $x$.  \qed

\begin{thm}
\label{thm:grtisn}
If $G$ is a graph of order $n$, then $\grt(G)=n$ if and only if there exists an integer $k$ such that
$n=2k$, and the vertices of $G$ can be labeled $x_1,\ldots,x_k,y_1,\ldots,y_k$ in such a way that \2 \\
\indent $\bullet$ $x_i$ is adjacent to $y_i$ for each $i$, \\
\indent $\bullet$ $\{x_1,\ldots,x_k\}$ is an independent set, and
\\
\indent $\bullet$ $y_j$ is adjacent to $x_i$ implies $i \ge j$.
\end{thm}
\proof  If the vertices of a graph $G$ can be labeled as in the statement of the theorem, then it is straightforward to check that $S=(x_1,\ldots,x_k,y_k,\ldots,y_1)$ is a total dominating sequence in $G$, and hence $\grt(G)=n$.

For the converse we assume that $\grt(G)=n$.  Let $S=(s_1,\ldots,s_n)$ be any Grundy total dominating sequence of length $n$ in $G$.  The first vertex of $S$ must have degree~$1$ since it footprints exactly one vertex.  Label this first vertex $x_1$ and label the vertex it footprints~$y_1$. By Lemma~\ref{lem:reciprocal} we know that $y_1$ footprints $x_1$.  Delete $x_1$ and $y_1$ from the sequence $S$ and label the first vertex that appears in the resulting sequence $x_2$.  Let $y_2$ be the unique vertex that $x_2$ footprints; that is,  $y_2=f_S^{-1}(x_2)$.  By Lemma~\ref{lem:reciprocal} $y_2 \not\in \{x_1,y_1\}$ and $y_2=f_S(x_2)$.  Once again we delete both $x_2$ and $y_2$ from the sequence.  Continuing in this fashion by choosing the first vertex of the remaining sequence to be $x_i$, denoting $f_S^{-1}(x_i)$ by $y_i$ and reasoning as above, we see that all  vertices of $G$ will be labeled and deleted from the sequence.  It follows that $n$ is even, say $n=2k$.  Let $X=\{x_1,\ldots,x_k\}$ and let $Y=\{y_1,\ldots,y_k\}$. By the way the vertices were labeled it is  clear that $x_i$ and $y_i$ are adjacent whenever $1\le i\le k$. Suppose $1 \le i < j \le k$. Since $x_i$ does not footprint $x_j$ with respect to $S$ we know that $x_i$ and $x_j$ are not adjacent in $G$.  Thus, $X$ is independent.  Moreover, $y_j$ is footprinted by $x_j$ so is not adjacent to $x_i$ for $i < j$.~\qed

\medskip
Clearly, the graphs from Theorem~\ref{thm:grtisn} all contain a perfect matching (it is given by the edges $x_iy_i$, for $i=1,\ldots,k$).  In an arbitrary graph the Grundy total domination number need not be bounded from below by the number of vertices in a matching.  The class of complete multipartite graphs contains graphs with  arbitrarily large matching number, and yet each has Grundy total domination number~$2$.  However, the Grundy total domination number is bounded below by the number of vertices in a semistrong matching. Indeed, the next result follows immediately by applying Observation~\ref{induced} and Theorem~\ref{thm:grtisn}.

\begin{cor} \label{cor:semistrong}
If $G$ is a graph, then $\grt(G) \ge 2 \nu_{ss}(G)$.
\end{cor}

In the next result we characterize the graphs with the smallest possible total Grundy domination number,  which is~$2$.

\begin{thm}
If $G$ is a graph, then $\grt(G)=2$ if and only if $G$ is a complete multipartite graph.
\label{thm:multi}
\end{thm}
\proof  It is clear that if $G$ is a complete multipartite graph, then $\grt(G)=2$.

For the converse, let $G$ be a graph with $\grt(G)=2$. We first note that if $x$ and $y$ are adjacent vertices, then $\{x,y\}$ is a (total) dominating set of $G$, since $S=(x,y)$ is a legal sequence, and so it must be a total dominating sequence. Next, we observe that if $x$ and $y$ are two nonadjacent vertices, then $N(x)=N(y)$. Indeed, otherwise $S=(x,y)$ or $S'=(y,x)$ would be a legal sequence but not total dominating sequence, because $x$ and $y$ are not totally dominated by the set $\{x,y\}$. Hence, if $A$ is a maximal independent set in $G$, then the neighborhoods $N(x)$ of all vertices $x$ from $A$ coincide.

Let $A$ be a maximal independent set in $G$, and $x\in A$. We claim that $N(x)\cup A=V(G)$. Suppose that there is a vertex $y\notin N(x)\cup A$. Since $y$ is not adjacent to $x$, it must be adjacent to some $x'\in A$, otherwise $A$ would not be a maximal independent set. But then $x$ and $x'$ are two nonadjacent vertices with $N(x)\ne N(x')$, a contradiction. Therefore, $N(x)\cup A=V(G)$.

Since $A$ was an arbitrarily chosen maximal independent set, we conclude that $G$ can be partitioned into maximal independent sets, each of which is adjacent to all other vertices not in that set. This implies that $G$ is the join of these maximal independent sets, and is thus a complete multipartite graph.~\qed

\section{Trees}
\label{sec:trees}

A \emph{rooted tree} distinguishes one vertex $r$ called the \emph{root}. Following the notation of~\cite{MHAYbookTD}, for each vertex $v \ne r$ of $T$, the \emph{parent} of $v$ is the neighbor of $v$ on the unique $(r,v)$-path, while a \emph{child} of $v$ is any other neighbor of $v$.
As observed earlier, the graphs $G$ satisfying $\grt(G) = n$ all contain a perfect matching. In the case of trees this condition is also sufficient.

\begin{thm}
If $T$ is a tree of order~$n$, then $\grt(T) = n$ if and only if $T$ has a perfect matching.
\label{cor:tree}
\end{thm}
\proof If $\grt(T) = n$, then by Theorem~\ref{thm:grtisn}, $T$ has a perfect macthing. Hence it suffices for us to show that for a tree $T$ with a perfect matching, $\grt(T) = n$. Let $T$ be such a tree with a perfect matching $M$, and note that $n=2k$ for some natural number $k$.
Choose an arbitrary vertex $r$ of $T$ and make $T$ a rooted tree with root $r$. Let $L$ be a linear order of
vertices of $T$ with the property that every child appears in the order before its parent (it is easy to see
that $L$ exists and also easy to construct it, by starting with leaves of $T$, deleting them,
and continuing in the same fashion). Now, let $(c_1,\ldots,c_k)$ be the suborder of $L$, obtained by choosing
the vertices from $T$ that are matched with respect to $M$ with their parent. Similarly, let $(p_k,\ldots,p_1)$
be the suborder of the dual $L^d$ of $L$ (obtained by reversing the order of $L$), in which the vertices
from $T$ that are matched with respect to $M$ with their child are chosen.
It is easy to see that $S=(c_1,\ldots,c_k,p_k,\ldots,p_1)$ is total dominating sequence of $T$
(each $c_i$ footprints its parent, and each $p_j$ footprints its child that it is matched to with respect to $M$).~\qed

\medskip
As an immediate consequence of Observation~\ref{induced} and Theorem~\ref{cor:tree}, we determine the Grundy total domination number of a path.

\begin{cor}
For $n \ge 2$ even, $\grt(P_n) = n$, while for $n \ge 3$ odd, $\grt(P_n) = n-1$.
\label{cor:path}
\end{cor}

We remark that a tree is open twin-free if and only if it has no strong support vertex. We next determine a lower bound on the Grundy total domination number of a tree with no strong support vertex in terms of its order, and we characterize the trees that achieve this lower bound. We remark that the requirement that the tree has no strong support vertex is essential here. For example, a star $K_{1,n}$ has Grundy total domination number~$2$ and therefore there is no constant $c > 0$ such that $\grt(T) \ge c|V(T)|$ for every star $T$.

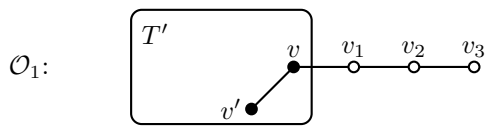
\begin{figure}[htb]
\begin{center}
\begin{tikzpicture}[scale=.8,style=thick,x=1cm,y=1cm]
\def\vr{2.5pt} 
\path (2.7,.75) coordinate (v);
\path (2.0,.05) coordinate (w);
\path (3.7,.75) coordinate (v1);
\path (4.7,.75) coordinate (v2);
\path (5.7,.75) coordinate (v3);
%
\draw (v) -- (w);
\draw (v) -- (v1);
\draw (v1) -- (v2);
\draw (v2) -- (v3);
\draw (v) [fill=black] circle (\vr);
\draw (w) [fill=black] circle (\vr);
\draw (v1) [fill=white] circle (\vr);
\draw (v2) [fill=white] circle (\vr);
\draw (v3) [fill=white] circle (\vr);
\draw [rounded corners] (0,-.2) rectangle (3,1.7);
\draw (.4,1.3) node {$T'$};
\draw (-1.7,.75) node {$\cO_1$:};
\draw[anchor = south] (v) node {$v$};
\draw[anchor = east] (w) node {$v'$};
\draw[anchor = south] (v1) node {$v_1$};
\draw[anchor = south] (v2) node {$v_2$};
\draw[anchor = south] (v3) node {$v_3$};
\end{tikzpicture}
\end{center}
\caption{The operation $\cO_1$.} \label{f:cO1}
\end{figure}

For this purpose, we define a family~$\cT$ of trees as follows. Let $\cT$ be the family of trees that contain a path $P_2$ and are closed under the operation $\cO_1$, which extends a tree $T'$ by adding a path $v_1v_2v_3$ and the edge $vv_1$ to a support vertex $v$ in the tree $T'$.
The operation $\cO_1$ is illustrated in Figure~\ref{f:cO1}, where here $v'$ is a leaf-neighbor of $v$ in $T'$. We remark that if $T \in \cT$ has order~$n$, then $n \equiv 2 \, (\modo \, 3)$. The two smallest trees in the family~$\cT$ are the path $P_2$ and the path $P_5$.

%

\begin{prop} \label{prop:cT}
If $T \in \cT$ has order~$n$, then $\grt(T) = \frac{2}{3}(n+1)$.
\end{prop}
\proof We proceed by induction on the order~$n \ge 2$ of a tree $T \in \cT$. If $n = 2$, then $T = P_2$ and $\grt(T) = 2 = \frac{2}{3}(n+1)$. This establishes the base case. Suppose that $n > 2$ and every tree $T' \in \cT$ of order~$n' < n$ satisfies $\grt(T') = \frac{2}{3}(n'+1)$. Let $T \in \cT$ be a tree of order~$n$. By construction of the family~$\cT$, the tree $T$ is obtained from a tree $T' \in \cT$ by adding a path $v_1v_2v_3$ and the edge $vv_1$ to a support vertex $v$ in the tree $T'$. Let $v'$ be a leaf-neighbor of $v$ in the tree $T'$. Let $T'$ have order~$n'$, and so $n' = n - 3$.

Every total dominating sequence of $T'$ can be extended to a total dominating sequence of $T$ by adding to it the vertices $v_1$ and $v_2$, and so $\grt(T) \ge \grt(T') + 2$. Conversely, suppose that $S$ is a Grundy total dominating sequence of $T$ (of maximum length). Suppose that the vertex $v_1$ appears in the sequence $S$. The vertex $v_1$ footprints $v$ or $v_2$. If $v_1$ footprints $v$, then the leaf $v'$ does not appear in the sequence $S$ and we can replace $v_1$ in $S$ with the vertex $v'$ (and leave all other entries in the sequence unchanged). If $v_1$ footprints $v_2$, then the leaf $v_3$ does not appear in the sequence $S$ and we can replace $v_1$ in $S$ with the vertex $v_3$ (and leave all other entries in the sequence unchanged). In both cases, we produce a new legal sequence, $S^*$, of $T$. If $S^*$ is not a total dominating sequence of $T$, then it can be extended to a total dominating sequence of $T$, contradicting the fact that $S$ is a Grundy total dominating sequence. Hence, $S^*$ is a total dominating sequence of $T$, implying that it is a Grundy total dominating sequence. Therefore, we can choose the sequence $S$ so that $v_1$ does not appear in the sequence $S$. With this choice of the sequence $S$, both $v_2$ and $v_3$ appear in the sequence $S$. Removing the vertices $v_2$ and $v_3$ from $S$ produces a total dominating sequence of $T'$, implying that $\grt(T') \ge \grt(T) - 2$. Consequently, $\grt(T) = \grt(T') + 2$. Applying the inductive hypothesis to the tree $T' \in \cT$, $\grt(T') = \frac{2}{3}(n'+1) = \frac{2}{3}(n-2)$, and so $\grt(T) = \grt(T') + 2 = \frac{2}{3}(n+1)$.~\qed

\medskip
Recall that if $T$ is a tree of order~$n \ge 3$, then $\gt(T) \le \frac{2}{3}n$. In contrast, we show next that $\grt(T) > \frac{2}{3}n$ for a tree $T$ of order~$n \ge 3$ with no strong support vertex.

\begin{thm} \label{thm:tree}
If $T$ is a nontrivial tree of order~$n$ with no strong support vertex, then
$\grt(T) \ge \frac{2}{3}(n+1)$,
with equality if and only if $T \in \cT$.
\end{thm}
\proof We proceed by induction on the order~$n \ge 2$ of a nontrivial tree. If $n = 2$, then $T = P_2$,  $\grt(T) = 2 = \frac{2}{3}(n+1)$ and $T \in \cT$. This establishes the base case. Assume that $n \ge 3$ and every nontrivial tree of order less than~$n$ with no strong support vertex satisfies the statement of the theorem. Let $T$ be a nontrivial tree of order~$n$ with no strong support vertex. Since $T$ has no strong support vertex, $T$ is not a star and $\diam(T) \ge 3$. If $\diam(T) = 3$, then $T = P_4$ and, by Corollary~\ref{cor:path}, $\grt(T) = 4 > \frac{2}{3}(n+1)$. Hence we may assume that $\diam(T) \ge 4$. In particular, $n \ge 5$.

We now root the tree $T$ at a leaf $r$ on a longest path in $T$. Let $u$ be a vertex at maximum distance from~$r$. Necessarily, $u$ is a leaf. Thus, $d_T(u,r) = \diam(T) \ge 4$. Let $v$ be the parent of $u$, let $w$ be the parent of $v$, and let $x$ be the parent of $w$. Since $u$ is a vertex at maximum distance from the root~$r$, every child of $v$ is a leaf. By supposition, $T$ has no strong support vertex, and so $d_T(v) = 2$.

Let $T_1 = T - \{u,v\}$ and let $T_1$ have order~$n_1$. Then, $n_1 = n - 2 \ge 3$. Suppose that the tree $T_1$ has no strong support vertex. If $n_1 = 3$, then $T_1 \cong P_3$ and $T_1$ has a strong support vertex, a contradiction. Hence, $n_1 \ge 4$. Applying the inductive hypothesis to $T_1$, $\grt(T_1) \ge \frac{2}{3}(n_1+1) = \frac{2}{3}(n-1)$. Let $S' = (v_1,\ldots,v_k)$ be a Grundy total dominating sequence of $T_1$, and so $k = \grt(T_1)$. Then the sequence $S = (u,v_1,\ldots,v_k,v)$ is a total dominating sequence of $T$ of length~$k + 2 = \grt(T_1) + 2 \ge \frac{2}{3}(n-1) + 2 = \frac{2}{3}(n+2)$, implying that $\grt(T) > \frac{2}{3}(n+1)$.  Hence, we may assume that the tree $T_1$ has a strong support vertex, for otherwise $\grt(T) > \frac{2}{3}(n+1)$, as desired.

Since $T$ has no strong support vertex but the tree $T_1$ has a strong support vertex, the vertex $w$ is necessarily a leaf in $T_1$ and its parent, namely $x$, is a strong support vertex in $T_1$. Let $w'$ be the leaf-neighbor of $x$ in $T_1$ different from $w$. We now consider the tree $T' = T - \{u,v,w\}$. Let $T'$ have order~$n'$, and so $n' = n - 3 \ge 2$. If $T'$ has a strong support vertex, then such a vertex is also a strong support vertex of $T$, a contradiction. Hence, $T'$ has no strong support vertex. Applying the inductive hypothesis to $T'$, $\grt(T') \ge \frac{2}{3}(n'+1) = \frac{2}{3}(n-2)$. Let $S' = (v_1,\ldots,v_k)$ be a Grundy total dominating sequence of $T'$, and so $k = \grt(T')$. Then, the sequence $S = (u,v_1,\ldots,v_k,v)$ is a total dominating sequence of $T$ of length~$k + 2 = \grt(T') + 2 \ge \frac{2}{3}(n-2) + 2 = \frac{2}{3}(n+1)$, implying that $\grt(T) \ge \frac{2}{3}(n+1)$. This establishes the desired lower bound.

Finally, suppose that $\grt(T) = \frac{2}{3}(n+1)$. This implies that $\grt(T') = \frac{2}{3}(n'+1)$. By the inductive hypothesis, $T' \in \cT$. Since $T$ can be obtained from the tree $T' \in \cT$ by applying operation~$\cO_1$ to the support vertex $x$ of $T'$, the tree $T \in \cT$. Conversely, by Proposition~\ref{prop:cT}, if $T \in \cT$ has order~$n$, then $\grt(T) = \frac{2}{3}(n+1)$.~\qed

\medskip
We remark that the result in Theorem~\ref{thm:tree} does not hold for general bipartite graphs that are open twin-free.
For $k \ge 2$, consider the bipartite graph $G_k$ formed by taking as one partite set, a set $A$ of $2k-1$ elements, and as the other partite set a set $B$ whose vertices correspond to all the $k$-element subsets of $A$, and joining each element of $A$ to those subsets to which it belongs. Let $S$ be a Grundy total dominating sequence of $G_k$. Every set of $k$ vertices chosen from $A$ totally dominates the set $B$, and so $S$ contains at most $k$ vertices of $A$. The first vertex of $B$ that appears in $S$ totally dominates~$k$ vertices of $A$. At most~$k-1$ additional vertices of $B$ appear in the sequence $S$ in order to totally dominate the remaining $k-1$ vertices of $A$. Therefore, the sequence $S$ contains at most~$k$ vertices from $B$, and so $\grt(G_k) \le 2k$. It is a simple exercise to show that $\grt(G_k) \ge 2k$, implying that $\grt(G_k) = 2k$. However, $G_k$ has order~$2k-1 + {2k-1 \choose k}$ and minimum degree~$\delta(G_k) = k$. This implies that no minimum degree is sufficient to guarantee that the Grundy total domination number of a general bipartite graph that is open twin-free is bounded below by a constant times its order. We state this formally as follows.

\begin{ob}
\label{ob:bip}
There is no constant $c > 0$ such that $\grt(G) \ge c|V(G)|$ for every bipartite graph $G$ that is open twin-free.
\end{ob}

\section{Regular Graphs}
\label{sec:regular}

In this section, we establish a lower bound on the Grundy total domination number of a $k$-regular graph. It is only of interest to consider values of $k \ge 2$ since if $k=1$, then $\gt(G) = \grt(G) = n$. We begin by determining the Grundy total domination number of a $2$-regular graph.

\begin{prop}
For $n \ge 3$ odd, $\grt(C_n) = n-1$, while for $n \ge 4$ even, $\grt(C_n) = n-2$.
\label{prop:cycle}
\end{prop}
\proof For $n \ge 3$, let $C_n$ be a cycle given by $v_1v_2 \ldots v_nv_1$. Suppose, firstly, that $n$ is odd. Then, $C_n$ contains as an induced subgraph a path $P_{n-1}$ on an even number of vertices. Thus, by Observation~\ref{induced} and Corollary~\ref{cor:path}, $\grt(C_n) \ge \grt(P_{n-1}) = n-1$. Conversely, since the first vertex in every total dominating sequence of $C_n$ footprints two vertices, we note that $\grt(C_n) \le n-1$. Consequently, $\grt(C_n) = n-1$ for $n$ odd. Suppose next that $n \ge 4$ is even. Let $A$ and $B$ be the two partite sets of $C_n$, and let $S$ be a Grundy total dominating sequence of $C_n$. The first vertex in $S$ that footprints a vertex of $A$ belongs to $B$ and footprints two vertices of $A$. The first vertex in $S$ that footprints a vertex of $B$ belongs to $A$ and footprints two vertices of $B$. Thus, at least two vertices in $S$ footprint two vertices, implying that $\grt(C_n) \le n-2$. Since $(v_1,v_2,\ldots,v_{n-2})$ is a total dominating sequence of $C_n$, $\grt(C_n) \ge n-2$. Consequently, $\grt(C_n) = n-2$ for $n$ even.~\qed

\medskip
As a consequence of a result due to Chv\'{a}tal and McDiarmid~\cite{ChMc}, if $G$ is a $k$-regular graph of order~$n$, then $\gt(G) \le (\lfloor
\frac{k+2}{2} \rfloor  / \lfloor \frac{3k}2 \rfloor) \, n$. In~\cite{SoHe-2013} it was shown that for every $k$-regular graph $G$ of order~$n$ with no isolates, $\Gamma_t(G) \le n/(2-\frac{1}{k})$. In contrast to these upper bounds on the total and upper total domination numbers of regular graphs, we establish a lower bound on the Grundy total domination number of a regular graph. As a special case of Proposition~\ref{prop:Delta}, we have the following lower bound on the Grundy total domination number of a regular graph.

\begin{cor}
For $k \ge 1$, if $G$ is a $k$-regular connected graph of order~$n$, then $\grt(G) \ge \frac{n}{k}$, with equality if and only if $G = K_{k,k}$.
\label{cor:reg}
\end{cor}

We show next that the lower bound of Corollary~\ref{cor:reg} can be improved considerably if the $k$-regular graph is different from $K_{k,k}$. By Proposition~\ref{prop:cycle}, it is only of interest to consider the case when $k \ge 3$.

\begin{thm}
For $k \ge 3$, if $G$ is a connected $k$-regular graph of order~$n$ different from $K_{k,k}$, then
\[
\grt(G) \ge \left\{
\begin{array}{ll}
\frac{n + \lceil \frac{k}{2} \rceil - 2}{k-1}  & \mbox{if $G$ is not bipartite} \3 \\
\frac{n + 2\lceil \frac{k}{2} \rceil - 4}{k-1}  & \mbox{if $G$ is bipartite.}
\end{array}
\right.
\]
\label{thm:regular}
\end{thm}
\proof For $k \ge 3$, let $G$ be a connected $k$-regular graph of order~$n$ different from $K_{k,k}$. We consider two cases.

\emph{Case~1. $G$ is not a bipartite graph.} Since $G$ is connected, there exists a pair of vertices of $G$ that are not open twins but have at least one common neighbor. Among all such pairs of vertices of $G$, let $v_1$ and $v_2$ be chosen to have the maximum number of common neighbors. We construct a total dominating sequence of $G$ as follows. Let $v_1$ and $v_2$ be the first and second vertices, respectively, in the sequence. We note that since $v_1$ and $v_2$ have at least one common neighbor, the vertex $v_2$ footprints at most~$k-1$ vertices. We now extend the subsequence $S_2 = (v_1,v_2)$ to a total dominating sequence of $G$ as follows.

Suppose that the current sequence is given by $S_i = (v_1,v_2,\ldots,v_i)$ for some $i \ge 2$. Let $B_i$ be the set of all vertices totally dominated by at least one vertex in $\widehat{S_i} = \{v_1,v_2,\ldots,v_i\}$. Suppose that $B_i \ne V(G)$, and so $\widehat{S_i} = \{v_1,v_2\ldots,v_i\}$ is not a total dominating set of $G$. We show that there must exist a vertex with at least one neighbor in $B_i$ and at least one neighbor not in $B_i$. Suppose this is not the case. Since $G$ is connected and $B_i \ne V(G)$, there is at least one vertex not in $B_i$ that is adjacent to a vertex of $B_i$. Let $A_i$ be the set of all vertices of $G$ that have a neighbor in $B_i$ but do not belong to the set $B_i$. By our supposition, every vertex in the set $A_i$ has all its $k$ neighbors in $B_i$. Further, $A_i \cap B_i = \emptyset$. If a vertex in $B_i$ has at least one neighbor in $A_i$ but fewer than $k$ neighbors in $A_i$, then such a vertex has a neighbor in $B_i$ and a neighbor not in $B_i$, a contradiction. Therefore, every vertex in $B_i$ that has a neighbor in $A_i$ has all its $k$ neighbors in $A_i$. The connectivity and the $k$-regularity of $G$ therefore imply that $G$ is a bipartite graph (with partite sets $A_i$ and $B_i$), a contradiction. Therefore, there is a vertex with at least one neighbor in $B_i$ and at least one neighbor not in $B_i$. As the $(i+1)$st vertex in our sequence, we choose such a vertex, say $v_{i+1}$, that footprints as few vertices as possible, and let $S_{i+1} = (v_1,\ldots,v_i,v_{i+1})$. Since $v_{i+1}$ has at least one neighbor in $B_i$, the vertex $v_{i+1}$  footprints at most~$k-1$ vertices.

Continuing in this way, let $S = S_t = (v_1,v_2,\ldots, v_t)$ be the final resulting sequence of length~$t$. Then, $S$ is a total dominating sequence. Further, the vertex $v_1$ footprints $k$ vertices, while every other vertex in the sequence footprints at most~$k-1$ vertices. We proceed further with the following claim.

\begin{unnumbered}{Claim~A.}
At least one of the vertices $v_2$ or $v_t$ footprints at most~$\lfloor \frac{k}{2} \rfloor$ vertices.
\end{unnumbered}
\textbf{Proof of Claim~A.} Suppose, to the contrary, that both $v_2$ and $v_t$ footprint at least~$\lfloor \frac{k}{2} \rfloor + 1$ vertices. In particular, since $v_2$ footprints at least~$\lfloor \frac{k}{2} \rfloor + 1$ vertices this implies that $v_1$ and $v_2$ have at most~$k - \lfloor \frac{k}{2} \rfloor - 1 = \lceil \frac{k}{2} \rceil - 1 \le \lfloor \frac{k}{2} \rfloor$ common neighbors.

We consider the final vertex in the sequence $S$, namely the vertex $v_t$. Let $U$ be the set of vertices footprinted by $v_t$ and let $F = V(G) \setminus U$. Thus, every vertex in $F$ is footprinted by some vertex $v_i$, where $1 \le i \le t-1$. By supposition, $|U| \ge \lfloor \frac{k}{2} \rfloor + 1$.

If a vertex $v \in U$ is adjacent to some other vertex of $U$, then we would have chosen the vertex $v$ instead of the vertex $v_t$ since $v$ footprints at most~$|U| - 1$ vertices which is fewer vertices than are footprinted by the vertex $v_t$, a contradiction. Therefore, $U$ is an independent set.

Let $v \in U$ and let $X = N(v)$. We note that $X \subseteq F$ and $|X| = k$. If a vertex $x \in X$ is not adjacent to every vertex in $U$, then we would have chosen the vertex $x$ instead of the vertex $v_t$ since $x$ footprints at most~$|U| - 1$ vertices once again contradicting our choice of the vertex $v_t$. Therefore, every vertex in $X$ is adjacent to every vertex in $U$.

Suppose that two vertices, $x_1$ and $x_2$ say, in $X$ are adjacent. Then, $x_1$ and $x_2$ are not open twins. Further, the vertices $x_1$ and $x_2$ have at least $|U| \ge \lfloor \frac{k}{2} \rfloor + 1$ common neighbors. However as observed earlier, $v_1$ and $v_2$ have at most~$\lfloor \frac{k}{2} \rfloor$ common neighbors. This contradicts our choice of $v_1$ and $v_2$ as a pair of vertices that are not open twins with the maximum number of common neighbors. Therefore, $X$ is an independent set.

Since $G \ne K_{k,k}$ and $X$ is an independent set, there must exist a pair of vertices in $X$ that are not open twins. However such a pair of vertices in $X$ have at least $|U| \ge \lfloor \frac{k}{2} \rfloor + 1$ common neighbors, once again contradicting our choice of $v_1$ and $v_2$. We deduce, therefore, that at least one of the vertices $v_2$ and $v_t$ footprints at most~$\lfloor \frac{k}{2} \rfloor$ vertices. This completes the proof of Claim~A.~\smallqed

\medskip
We now return to the proof of Case~1. As observed earlier, by the way in which the sequence $S$ is constructed, the vertex $v_1$ footprints $k$ vertices while every other vertex in the sequence footprints at most~$k-1$ vertices. By Claim~A, at least one of the vertices $v_2$ and $v_t$ footprints at most~$\lfloor \frac{k}{2} \rfloor$ vertices. We note that the number of footprinted vertices is precisely the order of the graph, namely~$n$. Thus, since the sequence $S$ has length~$t$, our earlier observations imply that $n \le k + \lfloor \frac{k}{2} \rfloor + (t-2)(k-1)$, or, equivalently, $\grt(G) \ge t \ge (n + \lceil \frac{k}{2} \rceil - 2)/(k-1)$. This completes Case~1.

\medskip
\emph{Case~2. $G$ is a bipartite graph.} Let $G$ have partite sets $A$ and $B$. We construct firstly a total dominating sequence, $S_A$, of $B$. Such a sequence $S_A$ satisfies $\widehat{S_A} \subseteq A$ and $\widehat{S_A}$ totally dominates the set $B$. Further, each vertex of the sequence $S_A$ footprints a vertex of $B$ not footprinted by any vertex preceding it in $S_A$.

Since $G$ is connected and $G \ne K_{k,k}$, there exists a pair of vertices in $A$ that are not open twins but have at least one common neighbor. Among all such pairs of vertices of $G$, let $a_1$ and $a_2$ be chosen to have the maximum number of common neighbors. We construct a total dominating sequence, $S_A$ of $B$ as follows. Let $a_1$ and $a_2$ be the first and second vertices, respectively, in the sequence. We note that since $a_1$ and $a_2$ have at least one common neighbor, the vertex $a_2$ footprints at most~$k-1$ vertices. We now extend the subsequence $(a_1,a_2)$ using the same selection method as in Case~1; that is, if the vertices in $B$ are not yet totally dominated by our chosen vertices selected from $A$, we choose the next vertex in the sequence so that it is adjacent to at least one vertex already footprinted and so that it footprints as few vertices (of $B$) as possible.
Let $S_A = (a_1,a_2,\ldots, a_t)$ be the resulting subsequence of $S$ such that every vertex in $B$ is footprinted by some vertex in $\widehat{S_A} = \{a_1,a_2,\ldots,a_t\}$.
An analogous, but slightly simpler proof to that presented in the proof of Claim~A in Case~1 shows that at least one of the vertices $a_2$ and $a_t$ footprints at most~$\lfloor \frac{k}{2} \rfloor$ vertices.

\begin{unnumbered}{Claim~B.}
At least one of the vertices $a_2$ or $a_t$ footprints at most~$\lfloor \frac{k}{2} \rfloor$ vertices.
\end{unnumbered}
\textbf{Proof of Claim~B.} Suppose, to the contrary, that both $a_2$ and $a_t$ footprint at least~$\lfloor \frac{k}{2} \rfloor + 1$ vertices. In particular, this implies that $a_1$ and $a_2$ have at most~$\lfloor \frac{k}{2} \rfloor$ common neighbors. Let $U$ be the set of vertices footprinted by $a_t$ and let $F = B \setminus U$. Thus, every vertex in $F$ is footprinted by some vertex $a_i$, where $1 \le i \le t-1$. By supposition, $|U| \ge \lfloor \frac{k}{2} \rfloor + 1$. Let $v \in U$ and let $X = N(v)$. We note that $X \subseteq A$ and $|X| = k$. If a vertex $x \in X$ is not adjacent to every vertex in $U$, then we would have chosen the vertex $x$ instead of the vertex $a_t$ since $x$ footprints fewer vertices than does $a_t$, contradicting our choice of the vertex $a_t$. Therefore, every vertex in $X$ is adjacent to every vertex in $U$. Since $G \ne K_{k,k}$, there must exist a pair of vertices in $X$ that are not open twins. However such a pair of vertices in $X$ have at least $|U| \ge \lfloor \frac{k}{2} \rfloor + 1$ common neighbors, contradicting our choice of the pair $a_1$ and~$a_2$.~\smallqed

\medskip
We now return to the proof of Case~2. By the way in which the sequence $S_A$ is constructed, the vertex $a_1$ footprints $k$ vertices while every other vertex in the sequence $S_A$ footprints at most~$k-1$ vertices. By Claim~B, at least one of the vertices $a_2$ and $a_t$ footprints at most~$\lfloor \frac{k}{2} \rfloor$ vertices. Analogously, we construct a total dominating sequence, $S_B = (b_1,b_2,\ldots,b_r)$, of $A$ consisting only of vertices of $B$ such that every vertex in $A$ is footprinted by some vertex in $\widehat{S_B}$. Further, the vertex $b_1$ footprints $k$ vertices while every other vertex in the sequence $S_B$ footprints at most~$k-1$ vertices and at least one of the vertices $b_2$ and $b_r$ footprints at most~$\lfloor \frac{k}{2} \rfloor$ vertices. Let $S$ be the sequence obtained by combining the sequences $S_A$ and $S_B$. Then, $S$ is a total dominating sequence of $G$. Exactly two vertices in the sequence $S$ footprint $k$ vertices, two vertices in $S$ footprint at most~$\lfloor \frac{k}{2} \rfloor$ vertices, and every other vertex in the sequence $S$ footprints at most~$k-1$ vertices. Therefore, since $S$ has length~$\ell = r + t$, $n \le 2k + 2(\lfloor \frac{k}{2} \rfloor) + (\ell - 4)(k-1)$, and so
\[
\grt(G) \ge \ell \ge \frac{n + 2\lceil \frac{k}{2} \rceil - 4}{k-1}.
\]
This completes the proof of Case~2, and of Theorem~\ref{thm:regular}.~\qed

\medskip
As an immediate consequence of Theorem~\ref{thm:regular}, we have the following result.

\begin{cor}
For $k \ge 3$, if $G$ is a connected $k$-regular graph of order~$n$ different from $K_{k,k}$, then $\grt(G) \ge \frac{n}{k-1}$, with strict inequality if $k \ge 5$.
\label{cor:regII}
\end{cor}

In the special case of Corollary~\ref{cor:regII} when $k = 3$ and $k = 4$, we have the following results. Recall that if $G$ is a cubic graph of order~$n$, then $\gt(G) \le \frac{1}{2}n$.

\begin{cor}
If $G \ne K_{3,3}$ is a connected cubic graph of order~$n$, then $\grt(G) \ge \frac{1}{2}n$. \label{cor:cubic}
\end{cor}

Recall that if $G$ is a $4$-regular graph of order~$n$, then $\gt(G) \le \frac{3}{7}n$.

\begin{cor}
If $G \ne K_{4,4}$ is a connected $4$-regular graph of order~$n$, then $\grt(G) \ge \frac{1}{3}n$. \label{cor:4reg}
\end{cor}

The connected cubic graph $G_4$ shown in Figure~\ref{f:prism}(b) of order~$n = 8$ satisfies $\grt(G_4) = 4 = \frac{1}{2}n$, while the connected $4$-regular graph $G_3$ shown in Figure~\ref{f:prism}(a) of order~$n = 12$ satisfies $\grt(G_3) = 4 = \frac{1}{3}n$. Thus, the bounds in Corollaries~\ref{cor:cubic} and~\ref{cor:4reg} are achievable.

\section{Relations between $\grt(G)$ and $\gr(G)$}
\label{sec:relations}

Let $S=(v_1,\ldots,v_k)$ be a sequence of distinct vertices of a graph $G$ and $\widehat{S}$ the corresponding set of vertices.
This sequence $S$ is called a \emph{legal (closed neighborhood) sequence}
if
\begin{equation}
\label{e:ordinary}
N[v_i] \setminus \bigcup_{j=1}^{i-1}N[v_j] \ne\emptyset,
\end{equation}
for every $i \in \{2,\ldots,k\}$.  If, in addition, $\widehat{S}$ is a dominating set of $G$,
then  $S$ is called a \emph{dominating sequence} of $G$.  The maximum length of a dominating sequence in $G$
is called the \emph{Grundy domination number} of  $G$ and is denoted by $\gr(G)$.  See~\cite{bgm-2014}.

As a direct analogy with the well-known inequality $\gamma_t(G)\le 2\gamma(G)$, which holds for an arbitrary
graph $G$ with no isolated vertices, we prove a general upper bound on $\grt(G)$ in terms of $\gr(G)$.

\begin{thm} \label{thm:interpolation}
If $G$ is a graph, then $\grt(G)\le 2\gr(G)$, and the bound is sharp.
\end{thm}
\proof  Let $S=(s_1,\ldots,s_k)$ be a total dominating sequence of $G$, where $k=\grt(G)$.
We will prove that at most $k/2$ vertices can be removed from $S$ in such a way that the resulting sequence $S'$ is
a legal closed neighborhood sequence of $G$.
We note that a vertex $s_i$ in the sequence $S$ prevents $S$ from being a legal closed neighborhood sequence  only if
$N[s_i] \setminus \cup_{j=1}^{i-1}N[s_j]=\emptyset$. Since $S$ is a total dominating sequence, we infer that in such a case
$s_i$ footprinted only vertices from $S$ that precede $s_i$. That is, $f_S^{-1}(s_i)\subseteq \{s_1,\ldots,s_{i-1}\}$,
where $f_S:V(G)\to \widehat{S}$ is a footprinter function, mapping each vertex to its footprinter.
Let
\[
A=\{s_i\in S \colon  f_S^{-1}(s_i)\subseteq \{s_1,\ldots,s_{i-1}\}\}.
\]

Suppose that for some vertex $s_j \in A$, the set $f_S^{-1}(s_j)\cap A$ is not empty. Let $s_i \in f_S^{-1}(s_j) \cap A$. Since $s_j \in A$, the vertex $s_i$ that is footprinted by $s_j$ satisfies $i < j$. Since $s_i \in A$, the vertex $s_i$ footprints some vertex $s_t$, where $t<i$. But then, when $s_j$ was added to $S$, $s_i$ was already totally dominated by $s_t$, a contradiction with $s_i \in f_S^{-1}(s_j)$. Therefore, for every vertex $s_j \in A$, the set $f_S^{-1}(s_j) \cap A = \emptyset$.

Suppose that $s_i,s_j \in A$, where $i<j$.  By definition of the set $A$, $f_S^{-1}(s_i)\subseteq \{s_1,\ldots,s_{i-1}\}$ and $f_S^{-1}(s_j)\subseteq \{s_1,\ldots,s_{j-1}\}$. Further, since every vertex is footprinted by a unique vertex in the sequence $S$, $f_S^{-1}(s_i)\cap f_S^{-1}(s_j)= \emptyset$. As observed earlier, the set $f_S^{-1}(s_j) \cap A = \emptyset$ for every $s_j\in A$. The collection of sets $\{f_S^{-1}(s_i) \colon s_i\in A\}$ therefore forms a partition of a subset of $\widehat{S}\setminus A$,
and for each $s_i\in A$, $|f_S^{-1}(s_i)|\ge 1$.  This readily implies that $|A|\le k/2$.  It now follows that the sequence $S'$,
which is obtained from $S$ by deleting vertices from $A$, is a legal closed neighborhood sequence of $G$. This sequence
$S'$ can be extended to a dominating sequence of $G$. Thus, $\gr(G)\ge k-|A| \ge k - \frac{k}{2} = \frac{k}{2} = \frac{1}{2}\grt(G)$. That the bound is sharp is shown by the class of complete graphs, $K_n$ with $n \ge 2$, that satisfy $\grt(K_n)=2$ and $\gr(G)=1$.~\qed

\medskip
On the other hand, a similar analog of the well-known lower bound $\gamma_t(G)\ge \gamma(G)$ does not hold for $\grt(G)$ in terms of $\gr(G)$.
Moreover, there exists no positive constant $c$ such that $\grt(G)\ge c\gr(G)$ would hold in general. For instance, if $G$ is the star $K_{1,n}$,
then $\grt(G)=2$, while $\gr(G)=n$.

\section{Edge covering sequences in hypergraphs}
\label{sec:hyper}

A connection between dominating sequences with covering sequences in hypergraphs
was established in \cite{bgm-2014}, and it can be done analogously for total dominating sequences.
We will use this connection to obtain two results, one about the possible lengths of total
dominating sequences, and the other about NP-completeness of the decision version of this problem.

Recall that given a hypergraph $\cH=(X,{\cal E})$ with no isolated vertices, an \emph{edge cover}
of $\cal H$ is a set of hyperedges from ${\cal E}$ that cover all vertices of $X$.  That is, the union of the hyperedges
from an edge cover is the ground set $X$. The minimum number of hyperedges in an edge cover
of $\cal H$ is called the \emph{(edge) covering number} of $\cal H$ and is denoted by $\rho(\cH)$,
cf.~\cite{BE}.
When a greedy algorithm is applied aiming to obtain an edge cover, hyperedges from $\cal H$
are picked one by one, resulting in a sequence ${\cal C}=(C_1,\ldots,C_r)$, where $C_i\in \cal E$.
In each step $i$, $1\le i \le r$,  $C_i$ is picked in such a way that it covers some vertex not captured by previous steps;
that is, $C_i\setminus (\cup_{j<i}{C_j})\neq \emptyset$.  If this condition is true for each
$i\in\{2,\ldots,r\}$, then we call ${\cal C}$ a \emph{legal (hyperedge) sequence} of $\cal H$.
If ${\cal C}=(C_1,\ldots,C_r)$ is a legal sequence and the set   $\widehat{\cal C}=\{C_1,\ldots,C_r\}$
is an edge cover of $\cal H$, then ${\cal C}$ is called an \emph{edge covering sequence}.
If the algorithm happens to produce an optimal solution, then $\widehat{\cal C}$
is a minimum edge cover of cardinality $\rho(\cal H)$, but in general $r\ge \rho(\cal H)$. The maximum length $r$ of an edge covering sequence of $\cal H$  is called the \emph{Grundy covering number} of $\cal H$, and
is denoted by $\rhogr(\cH)$.

Let $G$ be a graph with no isolated vertices, and let $\cH=(V(G),\cN(G))$, where $\cN(G)$ denotes the set of
all open neighborhoods of vertices in $G$, be the open neighborhood hypergraph of $G$. Clearly there is a
one-to-one correspondence between edge covering sequences  in the hypergraph $\cH$ and total dominating sequences
in $G$. Hence, using the following result from \cite{bgm-2014} we immediately
derive Corollary~\ref{cor:interpolation}.

\begin{thm} \label{thm:interpolation}
Let $\cH$ be a hypergraph.
For any number $\ell$ such that $\rho(\cH)\le \ell \le \rhogr(\cH)$ there
is an edge covering sequence of $\cal H$ having length $\ell$.
\end{thm}

\begin{cor} \label{cor:interpolation}
Let $G$ be a graph.
For any number $\ell$ such that $\gamma_t(G)\le \ell \le \grt(G)$ there
is a total dominating sequence of $G$ having length $\ell$.
\end{cor}


In the remainder of this section we will consider the following computational complexity problem:

\begin{center}
\fbox{\parbox{0.85\linewidth}{\noindent
{\sc Grundy Total Domination Number}\\[.8ex]
\begin{tabular*}{.93\textwidth}{rl}
{\em Input:} & A graph $G=(V,E)$, and an integer $k$.\\
{\em Question:} & Is $\grt(G)\ge k$?
\end{tabular*}
}}
\end{center}

It was shown in \cite{bgm-2014} that the {\sc Grundy Domination Number} problem,
which is the decision version of the Grundy domination number, is NP-complete,
by reduction from the following edge covering problem in hypergraphs

\begin{center}
\fbox{\parbox{0.85\linewidth}{\noindent
{\sc Grundy Covering Number in Hypergraphs}\\[.8ex]
\begin{tabular*}{.93\textwidth}{rl}
{\em Input:} & A hypergraph ${\cal H}=(X,\mathcal{E})$, and an integer $k$. \\
{\em Question:} & Is $\rhogr({\cal H})\ge k$?
\end{tabular*}
}}
\end{center}

\noindent that was also shown to be NP-complete.  (See~\cite[Theorem 4.2]{bgm-2014}.)
While the membership of {\sc Grundy Total Domination Number} in NP is trivial, we will show
the NP-hardness of this problem by a reduction from {\sc Grundy Covering Number in Hypergraphs}.

Given a hypergraph $\cH=(X,{\cal E})$  {\em the incidence graph} of $\cH$ is the bipartite graph $G=(V,E)$,
whose vertex set $V$ can be partitioned into disjoint (independent) sets $\widetilde{X}$ and $\widetilde{\cal E}$, which correspond to
the sets of vertices $X$ and hyperedges $\cal E$, respectively.  A vertex $\widetilde{x}\in \widetilde{X}$ is
adjacent to $\widetilde{B}\in \widetilde{\cal E}$ if and only if $x\in B$. Using the definitions, we easily see that $\widetilde{S}=(\widetilde{B_1},\ldots,\widetilde{B_t})$, a sequence of vertices from $\widetilde{\cal E}$,
is a total dominating  sequence of $\widetilde{X}$ if and only if $S=({B_1},\ldots,{B_t})$ is an edge covering sequence in $\cH$.
Hence the Grundy covering number of a hypergraph $\cH$ coincides with the maximum length of a total dominating sequence of
$\widetilde{X}$ in the incidence graph of $\cH$. In order to determine the Grundy total domination number of the incidence graph of $\cH$ we need to establish also the maximum length of a legal sequence of vertices from $\widetilde{X}$ that totally dominates $\widetilde{\cal E}$.  For this we introduce a new notion as follows.

Given a hypergraph $\cH=(X,{\cal E})$ a  {\em legal transversal sequence} is a sequence $S=(v_1,\ldots, v_t)$ of vertices from $X$ such that for each $i$ there exists an edge $E_i\in {\cal E}$ such that $v_i\in E_i$ and $v_j\notin E_i$ for all $j$, $j<i$.
That is, $v_i$ {\em hits} an edge $E_i$ which was not hit by any preceding vertices, and at the end every edge is hit by a vertex from $\widehat{S}$. The longest possible legal transversal sequence in a hypergraph $\cH$ will be called a {\em Grundy transversal sequence} and its length the {\em Grundy transversal number} of $\cH$, denoted $\taugr(\cH)$.

\begin{prop}
\label{p:rho-tau}
The Grundy transversal number of an arbitrary hypergraph $\cH$ equals the Grundy covering number of $\cH$; in symbols $$\taugr(\cH)=\rhogr(\cH).$$
\end{prop}
\proof
Let $S=(v_1,\ldots, v_t)$ be a Grundy transversal sequence in $\cH$, and let $(E_1,\ldots,E_t)$ be a legal sequence of edges in $\cH$ such that $E_i$ was hit by $v_i$ (i.e. $v_i\in E_i$) but was not
hit by the vertices that precede $v_i$ in $S$. We claim that the sequence $S'$ of these edges
in reverse order, that is $S'=(E_t,\ldots,E_1)$, is a Grundy covering sequence. Indeed, if $E_i$ is an arbitrary edge in this sequence, then $v_i\in E_i$, but $v_i\notin E_j$ for $j>i$, because in the transversal sequence $S$ the set $E_j$ was hit for the first time only later, by the vertex $v_j$. Hence, $S'$ is an edge covering sequence, and $\rhogr(\cH)\ge t=\taugr(\cH)$.

For the converse the same idea can be used. Note that if $S=(E_1,\ldots,E_u)$ is a Grundy covering sequence of $\cH$, and $v_i$ is a vertex that is in $E_i$
but is not in $E_1\cup\cdots\cup E_{i-1}$ (for each $i$), then the sequence $S'=(v_u,\ldots,v_1)$ is a legal transversal sequence. This implies $\taugr(\cH)\ge u=\rhogr(\cH)$.
\qed
\medskip

From Proposition~\ref{p:rho-tau}, and the fact that the Grundy total domination number of the incidence graph of $\cH$ coincides with $\taugr(\cH)+\rhogr(\cH)$ we derive the following result.

\begin{thm}
If $\cH$ is a hypergraph and $G$ the incidence graph of $\cH$, then $\grt(G)=2\rhogr(\cH)$.
\end{thm}

Since the reduction from a hypergraph to its incidence graph (which is bipartite) is efficiently computable,
it follows that {\sc Grundy Total Domination Number} is NP-hard even in bipartite graphs.

\begin{cor}
{\sc Grundy Total Domination Number} is NP-complete, even when restricted to bipartite graphs.
\end{cor}

\section{Open Problems}

We conclude with an open question and several open problems that we have yet to settle. By Observation~\ref{ob:bip}, there is no constant $c > 0$ such that $\grt(G) \ge c|V(G)|$ for every bipartite graph $G$ that is open twin-free. However, in our constructions every vertex belongs to a $4$-cycle. We pose the following question.

\begin{quest}
Does there exist a positive constant $c$ such that $\grt(G) \ge c|V(G)|$ for every bipartite graph $G$ with no $4$-cycles and with minimum degree at least~$2$?
\end{quest}

By Proposition~\ref{p:deg}, if $G$ is a graph of order~$n$, then $\grt(G) \le n - \delta(G) + 1$. We have yet to characterize the graphs achieving equality in this upper bound.

\begin{prob}
Characterize the graphs $G$ of order~$n$ for which $\grt(G) = n - \delta(G) + 1$.
\end{prob}

\begin{prob}
Find an efficient algorithm to compute the Grundy total domination number for trees.
\end{prob}

By Corollary~\ref{cor:cubic}, if $G \ne K_{3,3}$ is a connected cubic graph of order~$n$, then $\grt(G) \ge \frac{1}{2}n$. We observed that this bound is achievable.

\begin{prob}
Characterize the connected cubic graphs $G \ne K_{3,3}$ of order~$n$ for which $\grt(G) = \frac{1}{2}n$.
\end{prob}

By Corollary~\ref{cor:4reg}, if $G \ne K_{4,4}$ is a connected $4$-regular graph of order~$n$, then $\grt(G) \ge \frac{1}{3}n$. We observed that this bound is achievable.

\begin{prob}
Characterize the connected $4$-regular graphs $G \ne K_{4,4}$ of order~$n$ for which $\grt(G) = \frac{1}{3}n$.
\end{prob}

By Observation~\ref{chain}, for every graph $G$, $\gamma_t(G) \le \grt(G)$. By Theorem~\ref{thm:multi}, the graphs $G$ for which $\gt(G) = \grt(G) = 2$ are precisely the complete multipartite graphs. By Theorem~\ref{thm:Tdom3}, there is no graph $G$ satisfying $\gt(G) = \grt(G) = 3$. By Theorem~\ref{ob:gt4}, there are infinitely many connected graphs $G$ with $\grt(G) = \gt(G) = 4$. It remains an open problem to characterize the graphs $G$ for which $\gt(G) = \grt(G) = k$ when $k \ge 4$.

\begin{prob}
Characterize the graphs $G$ such that $\gt(G) = \grt(G) = k$ for $k \ge 4$.
\end{prob}

By Theorem~\ref{thm:interpolation}, if $G$ is a graph, then $\grt(G)\le 2\gr(G)$, and the class of complete graphs, $K_n$ with $n \ge 2$, achieve equality in this bound. However, it remains an open problem to characterize the extremal graphs.

\begin{prob}
Characterize the graphs $G$ for which $\grt(G)= 2\gr(G)$.
\end{prob}

\section{Acknowledgments}

The first author was supported by the Ministry of Science of Slovenia under the grant \\ P1-0297.  The
second author's research was supported in part by the South African National Research Foundation and
the University of Johannesburg. The research of the third author was supported by a grant from the Simons
Foundation (\#209654 to Douglas F. Rall).  This research was conducted during a visit to the University
of Maribor by the last two authors.



\begin{thebibliography}{99}

\bibitem{Alfewy} D. Archdeacon, J. Ellis-Monaghan, D. Fischer, D. Froncek, P.C.B. Lam, S. Seager, B. Wei, and R. Yuster. Some remarks on domination. \textit{J. Graph Theory} \textbf{46} (2004), 207--210.

\bibitem{BE} C.~Berge, \emph{Hypergraphs}, North-Holland, 1989, Amsterdam.

\bibitem{bgm-2014}
  B.~Bre\v{s}ar, T.~Gologranc, M.~Milani\v{c}, D.~F.~Rall, R.~Rizzi, Dominating sequences in graphs. \textit{Discrete Math.} \textbf{336} (2014), 22--36.

\bibitem{brklra-2010}
  B.~Bre{\v{s}}ar, S.~Klav{\v{z}}ar, and D.~F.~Rall,
  Domination game and an imagination strategy.
  \textit{SIAM J. Discrete Math.} \textbf{24} (2010), 979--991.

\bibitem{cs-79} C.~A.~Christen, S.~M.~Selkow, Some perfect coloring properties of graphs. \textit{J. Combin. Theory. Ser. B}  \textbf{27} (1979), 49--59.

\bibitem{ChMc} V. Chv\'{a}tal and C. McDiarmid, Small transversals in hypergraphs. \textit{Combinatorica} \textbf{12} (1992), 19--26.


\bibitem{CoDaHe80} E. J. Cockayne, R. M. Dawes, and S. T. Hedetniemi, Total domination in graphs. \textit{Networks} \textbf{10} (1980), 211--219.

\bibitem{dkr-15} P.~Dorbec, G.~Ko\v smrlj, and G.~Renault, The domination game played on unions of graphs.
 \textit{Discrete Math.}  \textbf{338} (2015),  71--79.

\bibitem{GH-2005} A.~Gy\'{a}rf\'{a}s and A.~Hubenko,  Semistrong edge coloring of graphs. \textit{J. Graph Theory}  \textbf{49}  (2005), 39--47.

\bibitem{hekara-2013}
  M.~A.~Henning, S. Klav\v zar, and D.~F.~Rall,
  Total version of the domination game.
  \textit{Graphs Combin.} \textbf{31}(5) (2015), 1453--1462.

\bibitem{hekara-2015}
  M.~A.~Henning, S. Klav\v zar, and D.~F.~Rall,
  The 4/5 upper bound on the game total domination number. \textit{Combinatorica}, to appear.

\bibitem{MHAYbookTD} M. A. Henning and A. Yeo, \emph{Total domination in graphs (Springer Monographs in Mathematics)}.  ISBN-13: 978-1461465249  (2013).

\bibitem{JRS-2014} F. Joos, D. Rautenbach, and T. Sasse, Induced matchings in subcubic graphs. \textit{SIAM J. Discrete Math.} \textbf{28}(1) (2014), 468--473.

\bibitem{KMM-2012}  R. J. Kang, M. Mnich, and T. M\"{u}ller, Induced matchings in subcubic planar graphs. \textit{SIAM J. Discrete Math.} \textbf{26} (2012), 1383--1411.


\bibitem{bill-2013}
  W.~B.~Kinnersley, D.~B.~West, and R.~Zamani,
  Extremal problems for game domination number.
  \textit{SIAM J. Discrete Math.} \textbf{27} (2013), 2090--2107.

\bibitem{SoHe-2013} J. Southey and M. A. Henning, Edge weighting functions on dominating sets. \textit{J. Graph Theory} \textbf{72} (2013), 346--360.

\bibitem{ThYe07} S. Thomass\'{e} and A. Yeo, Total domination of graphs and small transversals of hypergraphs. \textit{Combinatorica} \textbf{27} (2007), 473--487.

\bibitem{Tuza} Z. Tuza, Covering all cliques of a graph. \textit{Discrete Math.} \textbf{86} (1990), 117--126.


\end{thebibliography}
\end{document}